\documentclass[12pt,oneside]{article}

\usepackage{amsfonts}
\usepackage{amscd}
\usepackage{amsthm}
\usepackage{amssymb}
\usepackage{amsmath}
\usepackage{epsfig}

\title{Finite automorphisms of negatively curved Poincar\'e Duality groups}

\author{F. T. Farrell \footnote{This research was supported in part
by the National Science Foundation.} \hskip 5pt \& J.-F. Lafont}

\theoremstyle{definition}
\newtheorem{Def}{Definition}[section]
\newtheorem*{Rmk}{Remark}

\theoremstyle{proposition}

\newtheorem{Prop}{Proposition}[section]

\theoremstyle{plain}

\newtheorem{Thm}{Theorem}[section]
\newtheorem{Cor}{Corollary}[section]

\theoremstyle{remark}

\newtheorem*{Prf}{Proof}

\begin{document}

\maketitle

\begin{abstract}
\noindent In this paper, we show that if $G$ is a finite $p$-group ($p$ prime) acting by 
automorphisms 
on a $\delta$-hyperbolic Poincar\'e Duality group over $\mathbb Z$, then the fixed 
subgroup is a Poincar\'e Duality group over $\mathbb Z_p$.  We also provide a family of 
examples to show that the fixed subgroup might not be a Poincar\'e Duality group over 
$\mathbb Z$.  In fact, the fixed subgroups in our examples even fail to be duality groups
over $\mathbb Z$.

\end{abstract}

\section{Introduction.}

The study of finite group actions on topological spaces has a long and distinguished 
history.  A frequent theme is to try and understand the topology of the fixed point
set, both in it's intrinsic form, and as a subspace of the original space.  The classic
work of Smith shows that for finite $p$-groups acting on spheres, the fixed point set 
has the $\mathbb{Z}_p$ cohomology of a sphere.  However, there are examples of `exotic'
actions on spheres, where the fixed point set is not {\it homeomorphic} to a sphere
(indeed, does not even have the $\mathbb Z$ cohomology of a sphere).  
In this short paper, we are interested in relating actions on a hyperbolic group with the 
induced action on its boundary at infinity. 

We will start by relating the fixed subgroup of an automorphism with the fixed subset
of the induced action on the boundary at infinity.  In particular, this will allow us to
use the classic theorem of Smith to prove that if one starts with a $\delta$-hyperbolic 
Poincar\'e Duality 
group over $\mathbb Z$, and the group that is acting is a finite $p$-group ($p$ a prime),
then the fixed subgroup is a Poincar\'e Duality group over $\mathbb Z_p$.  We will 
then use the strict hyperbolization technique due to Charney and Davis [6] to construct 
examples of involutions of a Poincar\'e Duality group over $\mathbb Z$ whose fixed subgroup 
fails to be a Poincar\'e Duality group over $\mathbb Z$ (and in fact, aren't even duality
groups over $\mathbb Z$).  These examples also provide 
examples of `exotic' involution on a sphere (the boundary at infinity) which can be 
realized geometrically (i.e. by an isometry of a $CAT(-1)$ space).  They also show
that, in general, one could have involutions of $CAT(-1)$ spaces having a sphere 
as the boundary at infinity, where the induced involution on the boundary has a 
fixed point set which is {\it not} an ANR.

\vskip 10pt

\noindent {\bf Remark:} This paper was motivated by the following more specific
questions (each of which is still open).  Let $\Gamma =\pi_1(M)$ where $M$ is a
closed negatively curved Riemannian manifold, and let $\alpha:\Gamma \rightarrow
\Gamma$ be an automorphism with $\alpha ^2=Id_\Gamma$.

\vskip 5pt

\noindent {\bf Question 1}: Is the fixed subgroup $\Gamma ^\alpha$ a Poincar\'e 
Duality group over $\mathbb Z$?

\vskip 5pt

\noindent {\bf Question 2}: Is $\alpha$ induced by an involution of $M$?  That is to 
say, does there exist a self-homeomorphism $f:M\rightarrow M$ with $f^2=Id_M$, and
$f_\sharp =\alpha$?

\vskip 5pt

\noindent {\bf Question 3}: Let $\hat \alpha: \partial ^\infty \tilde M\rightarrow 
\partial
^\infty \tilde M$ be the induced involution of the sphere at infinity of the universal
cover of $M$.  Is the fixed point set of $\hat \alpha$ (when non-empty) an ANR?

\vskip 10pt

\centerline{\bf Acknowledgements.}

\vskip 5pt

\noindent The authors would like to thank the referee for several helpful comments, in 
particular for pointing out the existence of [5] and [8], and for suggesting the extension
of Theorem 2.2 that is included at the end of Section 2.  We would also like to thank the
members of the Geometric and Function Theory seminar at the University of Michigan for
pointing out a substantial simplification in our original proof of Proposition 2.1.

\section{Main results.}

\subsection{A positive result.}

\begin{Prop} \label{fixedsgpsbdry}
Let $\Gamma$ be a $\delta$-hyperbolic group, $\bar \sigma$ an automorphism of
$\Gamma$ of finite order $m$, and $\bar \sigma_\infty$ the induced action of $\bar \sigma$ on
$\partial ^\infty \Gamma$.  Then $(\partial ^\infty \Gamma)^{\bar \sigma_\infty}$ is
homeomorphic to $\partial ^\infty (\Gamma^{\bar \sigma})$.
\end{Prop}

\begin{Prf}
Let $\Sigma$ be a symmetric generating set for $\Gamma$, and consider the action of $\bar
\sigma$ on $\Gamma$.  Observe that if we define a new generating set $\Sigma^\prime :=
\bigcup_{i=1}^m \bar \sigma ^i(\Sigma)$, then $\bar \sigma$ acts by isometries on the Cayley
graph $Cay(\Gamma, \Sigma^\prime)$ of $\Gamma$ with respect to these generators.  Indeed, we
note that given any pair of elements $g, h$ in $\Gamma$, we have that:

$$d_{\Sigma^\prime}(g,h)=\inf \{i \hskip 5pt | \hskip 5pt g^{-1}h= \alpha_1\ldots \alpha_i,
\alpha _j\in \Sigma ^\prime \}$$

\noindent Taking a minimal such expression, and applying $\bar \sigma$ to it, we see that:

$$\bar \sigma(g)^{-1} \bar \sigma(h) = \bar \sigma (g^{-1}h)= \bar \sigma (\alpha_1) \ldots
\bar \sigma (\alpha _i)$$

\noindent But by invariance of $\Sigma ^\prime$ under $\bar \sigma$, we immediately get an
expression for $\bar \sigma(g)^{-1} \bar \sigma(h)$ as a product of $i$ elements of $\Sigma
^\prime$.  This forces $d_{\Sigma ^\prime}(\bar \sigma(g), \bar \sigma(h))\leq d_{\Sigma
^\prime}(g, h)$.  But now $\bar \sigma$, by hypothesis, has finite order $m$.  So by
iterating our inequality we get that:

$$d_{\Sigma ^\prime}(g, h) \geq d_{\Sigma ^\prime}(\bar \sigma(g), \bar \sigma(h)) \geq 
 \cdots \geq d_{\Sigma ^\prime}(\bar \sigma ^m(g), \bar \sigma ^m(h))= d_{\Sigma ^\prime}(g,
h)$$

\noindent which implies that all the inequalities are in fact equalities, and hence that
$\bar \sigma$ does indeed act by isometries on $Cay(\Gamma, \Sigma^\prime)$.  From now on, we
will omit the subscript $\Sigma ^\prime$ from our distance function in order to simplify
notation.

Our next step is to define certain subsets of the Cayley graph $Cay(\Gamma, \Sigma^\prime)$
in terms of their behavior under $\bar \sigma$, and to control the distance between these
subsets.  We let $F_k :=\{g \hskip 5pt | \hskip 5pt d(g, \bar \sigma g)\leq k\}$, and observe 
that $\Gamma ^{\bar \sigma} = F_0 \subseteq F_1 \subseteq F_2 \subseteq \cdots$.  Neumann 
[15] has shown that, for each i, there exists a $K_i$ such that $d(F_i, F_0)<K_i$.  

Next we observe that, by Neumann [15], the subgroup $\Gamma ^{\bar \sigma}=F_0$ is 
quasi-convex
in $\Gamma$.  In particular, $\partial ^\infty F_0$ embeds in $\partial ^\infty \Gamma$.
Now note that, trivially, we have that $\partial ^\infty (\Gamma ^{\bar \sigma})$ is in fact
a subset of $(\partial ^\infty \Gamma)^{\bar \sigma_\infty}$.  To prove equality, we need to
show the reverse inclusion.  So let us take a point $p\in (\partial ^\infty \Gamma)^{\bar 
\sigma_\infty}$, and let $\gamma \subset Cay(\Gamma, \Sigma^\prime)$ be a geodesic ray based
at the identity and with $\gamma (\infty)=p$.  Now by our choice of generators, we know that
$\eta := \bar \sigma (\gamma)$ will also be a geodesic ray (since $\bar \sigma$ acts
isometrically on the Cayley graph), and since the point $p=\gamma(\infty)$ is fixed by $\bar
\sigma_\infty$, we must have $d(\gamma, \eta)\leq C$ for some constant $C$. 

Our next claim is that, for each $n$, the inequality $d(\gamma (n), \eta (n))\leq 2C$ holds 
(and hence, as $\eta = \bar
\sigma (\gamma)$, forces $\gamma \subset F_{2C}$).  In order to see this, we consider
the following construction: given an integer $n$, we define $f(n)$ to be an integer
satisfying $d(\gamma (n), \eta (f(n)))\leq C$ (note that both $\gamma (n)$ and $\eta (f(n))$
correspond to elements in $\Gamma$).  We claim that $|f(n)-n|\leq C$ for all $n$.  By way
of contradiction, assume that $f(n)-n>C$.
The triangle inequality gives us:
$$f(n)=d(\eta(0), \eta(f(n)))\leq d(\eta(0), \gamma(n))+d(\gamma(n),\eta(f(n)))\leq n+C<f(n)$$
a contradiction (recall that $\eta(0),\gamma(0)$ are both the identity element in $\Gamma$).
The case $n-f(n)>C$ can be dealt with in an analogous manner.

We now know that, if $\gamma$ is an arbitrary geodesic ray originating at the identity, and
having $\gamma(\infty)=p$, then $\gamma
\subset F_{2C}$.  However, we also have that $d(F_0, F_{2C})\leq K$
for some constant $K$.  In particular, we can find a geodesic ray in $F_0$ which
has uniformly bounded distance from $\gamma$, which forces $p\in \partial ^\infty F_0 =
\partial ^\infty (\Gamma ^{\bar \sigma})$, completing the proof of the proposition.
\end{Prf}

\begin{Def}
We say that a topological space $X$ is an {$n$-dimensional Cech cohomology sphere with 
$R$ coefficients} (where $R$ is a PID) provided that $\check{\bar {H}} ^k(X; 
R)=0$ for all $k \neq n$, and $\check{\bar H} ^n(X; R)=R$ ($\check{ \bar
H}^*$ refers to reduced Cech cohomology).
\end{Def}

\begin{Def}
We say that a torsion-free group $G$ is a {\it duality group of dimension $n$ over $R$} 
(where again, $R$ is a PID), provided that there is a right $RG$-module $C$ such that
one has natural isomorphisms $H^k(G; A)\cong H_{n-k}(G; C\otimes _R A)$ for all 
$k\in \mathbb{Z}$ and all $RG$-modules $A$ (naturality is taken with respect to $A$, 
and $G$ acts diagonally on the tensor product $C\otimes _R A$).  If in addition we 
have that $C \cong R$, then we say that $G$ is a {\it Poincar\'e Duality group of 
dimension $n$ over $R$}.  Finally, if $G$ is a Poincar\'e Duality group of 
dimension $n$ over $R$, and the $G$ action on $C\cong R$ is trivial, we say that $G$ is
an {\it orientable Poincar\'e Duality group of dimension $n$ over $R$}.
\end{Def}

For background material on duality groups and Poincar\'e Duality groups, we refer to the 
lecture notes by Bieri [2].  Next, we quote the following result from Bestvina and Mess 
(Corollary 1.3 in their paper [1]):

\begin{Thm}[Bestvina \& Mess] \label{BEGM}
Let $\Gamma$ be a torsion-free $\delta$-hyperbolic group.  Then $\Gamma$ is a Poincar\'e 
Duality group of dimension $n$ over $\Lambda$ if and only if $\partial^\infty \Gamma$ 
is an $(n-1)$-dimensional Cech cohomology sphere with $\Lambda$ coefficients.
\end{Thm}

Using their result, we obtain an immediate corollary to our previous 
proposition:

\begin{Cor} \label{pdgrpsbdry}
Let $\Gamma$ be a torsion-free $\delta$-hyperbolic Poincar\'e Duality group of dimension 
$n$ over $\mathbb{Z}_p$. 
Let $G$ be a finite $p$-group ($p$ prime) acting by automorphisms on $\Gamma$.  Then there is 
a $0\leq k \leq n$ such that the subgroup $\Gamma ^G$ is a Poincar\'e Duality group of 
dimension $k$ over $\mathbb{Z}_p$.
\end{Cor}

\begin{Prf}
Let us first consider the case where $G$ is $\mathbb Z_p$.  Then
consider the induced action of $G$ on the boundary at infinity 
$\partial ^\infty \Gamma$.  Notice that, by Bestvina and Mess' result, $\partial ^\infty
\Gamma$ is a compact $(n-1)$-dimensional Cech cohomology sphere with $\mathbb{Z}_p$ 
coefficients.    
So we can use a version of Smith theory (see theorem III.7.11 in Bredon [4]), to get that the 
fixed point set of the action on the boundary at infinity must be a $(k-1)$-dimensional Cech 
cohomology sphere with $Z_p$ coefficients (for some $-1\leq k-1 \leq n-1$).  Now our 
Proposition \ref{fixedsgpsbdry} along with Bestvina and Mess' result immediately implies that 
the group $\Gamma ^G$ is a Poincar\'e Duality group of dimension $k$ over $\mathbb{Z}_p$ 
(where $0\leq k\leq n$).

For the more general case, we note that, since every $p$-group is solvable, one can find
a normal subgroup $G^\prime \leq G$. Finally one uses induction, since we have that $\Gamma 
^G = (\Gamma ^{G^\prime})^{G/G^\prime}$.  This gives the general case.
\end{Prf}

As was pointed out to the authors by the referee, Corollary 2.1 also follows from the
result announced by Chang and Skjelbred in [5], where they explain why the fixed set
of a finite $p$-group action on a Poincar\'e Duality space over $\mathbb Z_p$
is still a Poincar\'e Duality space over $\mathbb Z_p$. 

We conclude this section by mentioning that a Poincar\'e Duality group over $\mathbb Z$
is automatically a Poincar\'e Duality group over $\mathbb Z_p$, but that the converse
does not necessarily hold.  In Theorem 2.2, the group $\Gamma ^{\bar \sigma}$ will be
an example of this with $p=2$.

\subsection{A family of counterexamples.}

One can now ask the question of whether the previous result can be strengthened to obtain 
that the fixed subgroup is a Poincar\'e Duality group over $\mathbb{Z}$.  This turns out 
to be false, and in this section, we will construct counterexamples.  As was pointed out by
the referee, similar examples were constructed by Davis \& Leary [8].  Their 
construction used the reflection trick method (as opposed to our use of hyperbolization)
and served a somewhat different purpose.  We now proceed to state our main theorem.

\begin{Thm}
Let $\tau$ be a PL involution of a sphere $S^n$ whose fixed point set is a submanifold
$N^m$ which is {\it not} a homology sphere (with $\mathbb Z$ coefficients), and has dimension
$m\geq 2$.  Let $X$ be the 
strict hyperbolization of the 
suspension of $S^n$, and $\sigma$ the induced involution on $X$.  Let $\Gamma$ be the 
fundamental group of $X$, and $\bar \sigma$ the induced involution on $\Gamma$.  Then
$\bar \sigma$ is an involution of a ($\delta$-hyperbolic) orientable Poincar\'e Duality group 
over $\mathbb{Z}$ whose fixed subgroup $\Gamma ^{\bar \sigma}$ is
not a duality group over $\mathbb{Z}$.
\end{Thm}

Before starting with the proof, let us note that examples of involutions of spheres whose
fixed point sets are not homology spheres do exist.  In fact, Jones [11] has proved that every 
closed PL manifold that has the $\mathbb{Z}_2$ homology of a sphere can be realized as the 
fixed point set of a PL involution of some larger dimensional sphere.  

For a more concrete example, we can consider Brieskorn spheres: for $n\geq 2$, define two 
complex functions
$f_n(z_0,\ldots ,z_{2n+1}):=z_0^3+ \sum_{i=1}^{2n+1} z_i^2$, and $g_n(z_0,\ldots ,z_{2n}):= 
z_0^3+ \sum_{i=1}^{2n} z_i^2$.  Using these two functions, define a pair of manifolds $M_n$ 
and $N_n$ by considering
the intersection of $f_n^{-1}(0)$ and $g_n^{-1}(0)$ with a small enough ball centered at the 
origin in the 
appropriate complex vector space.  It is known that $M_n$ is PL homeomorphic to the sphere 
$S^{4n+1}$, while $N_n$ is a $(4n-1)$-dimensional manifold that does not have the 
$\mathbb{Z}$-homology of a sphere (combine Lemma 8.1 with the comments on pg. 72 in
Milnor [14]).  Furthermore, observe that the involution $z_{2n+1} 
\leftrightarrow -z_{2n+1}$ on $M_n$ has fixed point set $N_n$, giving us an infinite family
of examples.

\begin{Prf}
We start by recalling that the strict hyperbolization procedure given by Charney and Davis
(section 7 in [6])
takes a simplicial complex and functorially assigns to it a topological space (in fact,
a union of compact hyperbolic manifolds with corners)
that supports a metric of strict negative curvature.  Let us apply this procedure to the
suspension of the sphere $\Sigma S^n$ (respectively $\Sigma N^m$), and call the resulting 
space $X^{n+1}$ (respectively $Y^{m+1}$).  We will omit the dimension of the spaces unless
we explicitly require them for computations.

We now list out some properties of the spaces $X$ and $Y$. Observe that, by a result of 
Illman [9], there exists a triangulation of the pair $(S^n, N^m)$ such that the involution 
$\tau$ is a simplicial map.  In particular, the involution on the suspension will still be 
simplicial, and $\Sigma N^m$ is a subcomplex of $\Sigma S^n$.  Functoriality of 
the strict hyperbolization procedure now implies that $Y$ is
a totally geodesic subspace of $X$, invariant under the induced involution $\sigma$ on $X$ .  
Since hyperbolization preserves the local structure,
$X$ will be an orientable manifold, while $Y$ will have a pair of non-manifold points 
(corresponding to the two vertices of the suspension).  

Now take a basepoint $*\in Y\subset X$, and let $\Lambda =\pi_1(Y,*)$, $\Gamma =\pi_1(X,*)$. 
The involution $\sigma$ will give an order two automorphism $\bar \sigma$ of the group 
$\Gamma$.  We note that, since $\Gamma$ is the fundamental group of a closed orientable 
aspherical manifold, it is automatically an orientable Poincar\'e Duality group over 
$\mathbb{Z}$.  Now consider the fixed subgroup $\Gamma ^{\bar \sigma}$.  In order to get 
information about this group, we consider a lift of the action to the universal cover 
$\tilde X$ of $X$.

Let $\tilde *\in \tilde X$ be a preimage of the point $*$, and let us lift the involution 
$\sigma$
to the universal cover.  Note that the fixed point set of the lifted involution is precisely
the path connected lift $\tilde Y$ of $Y$ that contains the point $\tilde *$.  Furthermore, 
the
action $\bar \sigma$ on $\Gamma$ is compatible with the lift $\tilde \sigma$ of $\sigma$, 
in the sense that $(\bar \sigma (g))(\tilde *)=\tilde \sigma (g(\tilde *))$.

Next we note that $\Gamma ^{\bar \sigma} = \Lambda$.  Indeed $\Lambda$ is automatically fixed 
by $\bar \sigma$, hence we have a containment $\Lambda \leq \Gamma ^{\bar \sigma}$.  On the
other hand, for an arbitrary $g\in \Gamma ^{\bar \sigma}$, we have that $\tilde \sigma 
(g(\tilde *))=(\bar \sigma (g))(\tilde *)= g(\tilde *)$.  In particular, $g(\tilde *)$ must
be fixed under $\tilde \sigma$, which implies $g(\tilde *)\in \tilde Y$.  Since $\tilde Y$
is a path connected, totally geodesic subset, we can connect $\tilde *$ to $g(\tilde *)$
by a path which lies entirely within $\tilde Y$.  Looking at the projection of this path in
$X$, we observe that it is a closed loop based at $*$, representing the element $g$, and
lying entirely in $Y$.  Hence $g\in \Lambda$, giving us the reverse containment.  We conclude
that the two groups are equal.

\begin{figure}
\label{graph}
\begin{center}
\includegraphics[width=2in]{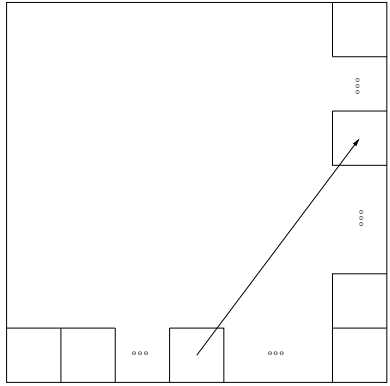}
\caption{Zeeman spectral sequence for our spaces.}
\end{center}
\end{figure}

So in particular, $Y$ is a topological space which happens to be a $K (\Gamma ^{\bar 
\sigma}, 1)$.  In particular, the group cohomology of $\Gamma ^{\bar \sigma}$ is related
to the compactly supported cohomology of $\tilde Y$.  So we have now reduced our claim to
analyzing the properties of $H^*(\tilde Y, \mathbb{Z})$.  In order to do this, we consider 
the Zeeman spectral sequence; let us first introduce some
terminology.  We will denote by $h_p$ the $p^{th}$ local homology sheaf for $Y$, and by 
$\tilde h_p$ the corresponding sheaf for $\tilde Y$.  For $x\in Y$ (respectively, in $\tilde 
Y$), we will denote by $h_p(x)$ (respectively $\tilde h_p(x)$) the stalk at the point $x$. 
Recall that $\tilde Y$ is the hyperbolization of an  
$(m+1)$-dimensional complex $\Sigma N^{m}$; we will use $Y^i$ to denote the subspace of 
$Y$ obtained from the hyperbolization of the $i$-skeleton of $\Sigma N^m$.  
Observe the following facts about the local homology sheaf:
\begin{itemize}
\item if $i\neq m+1$ and $x\notin Y^0$, then $h_i(x)=0$.
\item if $x \notin Y^0$, then $h_{m+1}(x)= \mathbb{Z}$.
\item there exists a point $p\in Y^0$ and an integer $s$ such that $2\leq s<m+1$ and $h_s(p)\neq 
0.$ 
\end{itemize}

\noindent  All of the previous remarks are clear, with the possible exception of the 
third: let $p$ be one of the two vertices of the suspension.  Since
the original link of $p$ was {\it not} a homology sphere, and as hyperbolization does
not change the link, there must exist an $s< m+1$ which yields the desired fact.
Note that the sheafs we are considering are given by local data, so that we have that
$\tilde h_p(\tilde x)=h_p(x)$, whenever $\tilde x$ is a lift of the point $x$.  Hence we
have that the three facts mentioned above for the stalks of the local homology 
sheaf $h_p$ on $Y$ also hold for the stalks of the local homology sheaf $\tilde h_p$
on $\tilde Y$.

Now the Zeeman spectral sequence (see section 2 of McCrory [13], based on previous work of 
Zeeman [17]) states that:
$$E^2_{i,j} := \check H_c^j(\tilde Y; \tilde h_i) \stackrel{j}{\Longrightarrow} 
H_{i-j}(\tilde Y)$$
with differentials $d^t_{i,j}:E_{i,j}^t\longrightarrow E^t_{i+(t-1),j+t}$.  Observe that, by 
the properties listed above for the $i^{th}$ local homology sheaf, $E_{i,j}^2=0$ if $i\neq 
m+1$ and $j\neq 0$.  So in particular, all the terms vanish except those in the $0^{th}$ 
row and those in the $(m+1)^{st}$ column (see figure on previous page).  

We now plan on working with this spectral sequence.  Observe from the shape of the spectral
sequence that one has isomorphisms:
$$E^2_{s,0}\cong E^3_{s,0} \cong \cdots \cong E^{m-s+2}_{s,0}$$
$$E^2_{m+1,m-s+2} \cong E^3_{m+1,m-s+2} \cong \cdots \cong E^{m-s+2}_{m+1,m-s+2}$$
and that the differential $d^{m-s+2}$ maps $E^{m-s+2}_{s,0}$ to $E^{m-s+2}_{m+1,m-s+2}$.
However, we know that $H_s(\tilde Y)=0$, so the differential must be an isomorphism.  This
yields:
$$H^0_c(\tilde Y; \tilde h_s) = E^2_{0,s} \cong E^2_{m+1,m-s+2} = H^{m-s+2}_c(\tilde Y; 
\tilde h_{m+1})$$
(since we are dealing with complexes, Cech cohomology coincides with standard cohomology).
Furthermore, $\tilde Y$ is simply-connected and has dimension $m+1\geq 3$, hence $\tilde 
h_{m+1}$ is the trivial $\mathbb Z$ sheaf over $\tilde Y^{m+1}-\tilde Y^0$.  This implies:
$$H^0_c(\tilde Y; \tilde h_s)\cong H^{m-s+2}_c(\tilde Y; \tilde h_{m+1}) = H^{m-s+2}_c(\tilde 
Y; \mathbb Z)$$
Now focusing on the left hand term, we note that $\tilde h_s(q)=0$ for all $q\notin \tilde Y
^0$, which gives us:
$$H^{m-s+2}_c(\tilde Y; \mathbb Z) \cong H^0_c(\tilde Y; \tilde h_s) = \bigoplus _{q\in 
\tilde Y^0} \tilde h_s(q)$$
But now observe that if $\tilde p$ is a vertex in $\tilde Y^0$ which is a lift of $p$ (one
of the vertex points of the suspension), then $\tilde h_s(g \cdot \tilde p)= h_s(p)
\neq 0$ for every element $g \in \Gamma ^{\bar \sigma}$.  Since all the points $g\cdot \tilde 
p$ lie in 
$\tilde Y^0$, and since $\Gamma ^{\bar \sigma}$ is an infinite group, this implies that 
$\bigoplus _{q\in 
\tilde Y^0} \tilde h_s(q)$ is not finitely generated.  So in particular, $H^{m-s+2}_c(\tilde 
Y; \mathbb Z)$ is not finitely generated.  Since $Y$ is a finite complex which happens to 
be a $K(\Gamma ^{\bar \sigma}, 1)$, we conclude that $H^{m-s+2}(\Gamma^{\bar \sigma}, \mathbb 
Z \Gamma^{\bar \sigma}) \cong H^{m-s+2}_c(\tilde Y; \mathbb Z)$ is
not finitely generated.  By Bieri and Eckmann's criterion (see  Bieri [2], section 9.10), 
this implies that $\Gamma ^{\bar
\sigma}$ cannot be a Poincar\'e Duality group over $\mathbb {Z}$.

In order to see that $\Gamma ^{\bar \sigma}$ is not even a duality group over $\mathbb Z$,
it is sufficient to show that the cohomological dimension of 
$\Gamma ^{\bar \sigma}$ is greater than $m-s+2$.  We first note that, 
since $s\geq 2$, we have that
$m-s+2\leq m$, so it is sufficient to show that $\Gamma ^{\bar \sigma}$ has non-trivial 
cohomology in some dimension
that is strictly greater than $m$.  Observe that, by construction, 
we have that $\Gamma ^{\bar \sigma}$ is the fundamental group of the finite aspherical 
$(m+1)$-dimensional space $Y$, which implies that the cohomological dimension of 
$\Gamma ^{\bar \sigma}$ is at most $m+1$.  We would be done provided we can show that the 
cohomological dimension of $\Gamma ^{\bar \sigma}$ is exactly $m+1$.  Looking back at the 
construction of $Y$, we observe that the submanifold $N^m$ we started with is a ${\mathbb Z}
_2$ homology sphere.  Suspending the manifold, we obtain an $(m+1)$-dimensional space
which is a ${\mathbb Z}_2$ homology manifold.  Now $Y$ is the hyperbolization of this space,
and since the hyperbolization procedure preserves the local structure, $Y$ is also an 
$(m+1)$-dimensional ${\mathbb Z}_2$ homology manifold.  This implies that 
$H^{m+1}(\Gamma ^{\bar \sigma}; \mathbb{Z}_2)\cong {\mathbb Z}_2 \neq 0$, which forces 
the cohomological dimension of $\Gamma ^{\bar \sigma}$ to be at least $m+1$.  This 
completes our proof.
\end{Prf}

\begin{Rmk}  As was pointed out to the authors by the referee, the argument in Theorem 2.2 
can also be used to
show that the condition that $G$ be a p-group in Corollary 2.1 really is necessary.  Namely,
there are examples of a $\mathbb Z_6$ action on an orientable Poincar\'e Duality 
$\delta$-hyperbolic group over $\mathbb Z$ whose fixed subgroup is not a Duality group
over {\it any} PID (in which $0\neq 1$).  Indeed, note that the unit tangent bundle 
$S(S^{n-1})$ 
of an $(n-1)$-dimensional sphere can be identified with the Stiefel manifold $V_{2,n}$ of 
orthonormal $2$-frames in $\mathbb R^n$.  The latter can be embedded in $\mathbb C^n$ 
via the map $f(u,v)=u+iv$ (where $u,v\in \mathbb R^n$ are orthonormal vectors).  

Note that since $u,v$ are orthonormal, we have that $|f(u,v)|^2=|u|^2+|v|^2=2$,
and also that:
$$\sum _{j=1}^n (u_j+iv_j)^2=\sum _{j=1}^n u_j^2 -\sum _{j=1}^n v_j^2 +2i(\sum _{j=1}^n
u_jv_v)=0$$
This implies that $S(S^{n-1})$ is diffeomorphic to the Brieskorn variety for the 
polynomial $z_1^2+\cdots z_n^2=0$.  In particular, we see that $S(S^{n-1})$ is the
fixed point set of the $\mathbb Z_6$ action on the Brieskorn variety for the polynomial
$z_0^3+z_1^2+\cdots z_n^2 + z_{n+1}^2=0$, where the action is given by $g(z_0,z_1,\ldots 
,z_n,z_{n+1})= (\theta z_0, z_1,\ldots ,z_n, -z_{n+1})$, where $\theta =e^{
2\pi i/3}$.  Furthermore, since odd dimensional spheres have a non-zero vector field, 
we have that $H_*(S(S^{2n-1}); \mathbb Z)\cong H_*(S^{2n-1}\times S^{2n-2}; \mathbb Z)$.  

Now let $N_n$ be the fixed point set of the above mentioned action of $\mathbb Z_6$ on
the Brieskorn variety $M_n$ for the polynomial $z_0^3+z_1^2+\cdots z_{2n}^2 + z_{2n+1}^2=0$.
As we mentioned earlier, the Brieskorn variety $M_n$ is PL-homeomorphic to $S^{4n+1}$,
while by the previous paragraph, $N_n$ is diffeomorphic to $S(S^{2n-1})$.  Suspending
the spaces and hyperbolizing gives us a $\mathbb Z_6$ action on a $CAT(-1)$ space $X$,
where now the fixed subset $Y$ is the hyperbolization of the suspension of $S(S^{2n-1})$.

The proof that $\pi_1(Y)$ is not Poincar\'e Duality over any PID $R$ is almost a verbatim
repetition of that given for Theorem 2.2.  In particular, the local homology sheaf for 
the space $Y$ will have three distinct 
indices (namely $s=2n-1,2n, 4n-2$) for which $h_s(p)\neq 0$ (where again, $p$ is one of
the suspension points).  Working through the Zeeman spectral sequence, we again find 
indices ($<4n-2$) where the cohomology of $\pi_1(Y)$ is not finitely generated.  The 
only substantial
change is in the argument showing that the cohomological dimension of $\pi_1(Y)$ over
$R$ is $4n-2$.  To do this, we merely note that the hyperbolization map $Y\rightarrow 
\Sigma S(S^{2n-1})$ induces a surjection on integral homology, together with the fact
that $H_{4n-2}(\Sigma S(S^{2n-1});R)\cong R$.
\end{Rmk}

\section{Concluding remarks.}

We finish our paper with a few remarks.  Firstly, 
we note that the results we obtain are, in some sense, dealing with exceptional 
automorphisms
of $\delta$-hyperbolic groups.  Indeed, Levitt \& Lustig [12] have shown that, in a suitable
sense, `most' automorphisms of a $\delta$-hyperbolic group have very simple fixed point 
sets for their induced actions on the boundary at infinity (in fact, their fixed point
sets consist of a pair of points).  Also, if we start with a torsion free group, then
the group of inner automorphisms will also be torsion free, hence any automorphism 
of finite order in some sense `lives' in the outer automorphism group, which tends to 
be small.

Secondly, we should point out that, in the counterexamples we constructed, the groups 
$\Gamma$ all have boundary at infinity which is in fact {\it homeomorphic} to a sphere.  
This follows from the fact that the link of every vertex in the space $\tilde X$ is 
PL-homeomorphic to the standard sphere, so by a result of Davis \& Januszkiewicz [7], the
boundary at infinity of $\tilde X$ is homeomorphic to a sphere.    

Thirdly, we can ask related questions in a somewhat more general setting.  More precisely, 
given an arbitrary topological space $Y$, we can consider the question of what type of 
actions can be realized {\it algebraically} or {\it geometrically}.  By a geometric action, 
we mean one that is induced by an isometry of a $\delta$-hyperbolic space $X$ whose boundary 
is homeomorphic to $Y$.  By an algebraic action, we mean one that is induced by an 
automorphism of a $\delta$-hyperbolic group $\Gamma$ whose boundary is homeomorphic to $Y$.
Note that, at the cost of changing the set of generators for the group $\Gamma$ (as in 
the proof of Proposition 2.1), we can always view an algebraic action as a geometric one 
(given by an isometry of the Cayley graph).

The fact that this question is non-trivial, even in the more general setting, can be seen
by considering the situation of a Menger manifold.  It is well known that there are 
numerous $\delta$-hyperbolic groups whose boundary at infinity are Menger manifolds.  
Now a result of Iwamoto [10] states that every closed subset of a Menger manifold can be 
realized as the fixed point set of an involution.  On the other hand, if an involution
can be realized algebraically via an involution $\sigma$ of a group $\Gamma$, then the fixed 
point set on the boundary at infinity must coincide with the boundary at infinity of 
the subgroup $\Gamma ^\sigma$.  However, the latter set cannot have any cutpoints (see 
Bowditch [3] and Swarup [1]).  This gives a necessary condition for a closed subset of a 
Menger manifold to be the fixed point set of an algebraically realizable involution.  What 
are the sufficient conditions?  

Finally, we mention that these examples give involutions of a $\delta$-hyperbolic group
$\Gamma$ where the fixed point set of the induced involution on the boundary at infinity
is not an ANR, although $\partial ^\infty \Gamma \cong S^n$.  One could ask whether the
fixed point set could display other complicated behavior.  For instance, does there 
exist an involution of a $\delta$-hyperbolic group $\Gamma$, with fixed subgroup 
$\Lambda$, with the property that $\partial ^\infty \Gamma \cong S^n$, $\partial ^\infty 
\Lambda \cong S^{n-2}$, and the embedding $S^{n-2}\cong \partial ^\infty \Lambda 
\hookrightarrow \partial ^\infty \Gamma \cong S^n$ is a locally flat, non-trivial knot?

\section{Bibliography.}

\vskip 10pt

\noindent [1] Bestvina, M. \&  Mess, G.  {\it The boundary of negatively curved groups}.  J.
Amer. Math. Soc.  4  (1991),  no. 3, pp. 469--481.

\vskip 5pt

\noindent [2] Bieri, R.  {\it Homological dimension of discrete groups}. Second edition.
Queen Mary College Mathematical Notes. Queen Mary College, Department of Pure Mathematics,
London,  1981. iv+198 pp.

\vskip 5pt

\noindent [3] Bowditch, B. H.  {\it Connectedness properties of limit sets}.
 Trans. Amer. Math. Soc.  351  (1999),  no. 9, pp. 3673--3686.

\vskip 5pt

\noindent [4] Bredon, G. E.  {\it Introduction to compact transformation groups}.
Pure and Applied Mathematics, Vol. 46.
Academic Press, New York-London,  1972. xiii+459 pp.

\vskip 5pt

\noindent [5] Chang, T. \& Skjelbred, T.  {\it Group actions on Poincar\'e Duality spaces}.
Bull. Amer. Math. Soc. 78 (1972), pp. 1024--1026.

\vskip 5pt

\noindent [6] Charney, R. M. \& Davis, M. W.  {\it Strict hyperbolization}.
 Topology  34  (1995),  no. 2, pp. 329--350.

\vskip 5pt

\noindent [7] Davis, M. W. \&  Januszkiewicz, T. {\it Hyperbolization of polyhedra}.
 J. Differential Geom.  34  (1991),  no. 2, pp. 347--388.

\vskip 5pt

\noindent [8] Davis, M. W. \& Leary, I. J.  {\it Some examples of discrete group actions
on aspherical manifolds}.  To appear in  The Proceedings of ICTP Conference on 
High-dimensional Manifold Topology.

\vskip 5pt

\noindent [9] Illman, S. {\it Smooth equivariant triangulations of $G$-manifolds for $G$ a
finite group}.  Math. Ann.  233  (1978),  no. 3, pp. 199--220.

\vskip 5pt

\noindent [10] Iwamoto, Y. {\it Fixed point sets of transformation groups of Menger manifolds,
their pseudo-interiors and their pseudo-boundaries}.
 Topology Appl.  68  (1996),  no. 3, pp. 267--283.

\vskip 5pt

\noindent [11] Jones, L.  {\it Construction of surgery problems}. in 
 Geometric topology (Proc. Georgia Topology Conf., Athens, Ga., 1977), 
 pp. 367--391, Academic Press, New York-London, 1979.

\vskip 5pt

\noindent [12] Levitt, G. \&  Lustig, M.  {\it Most automorphisms of a hyperbolic group have
very simple dynamics}.
 Ann. Sci. \'Ecole Norm. Sup. (4)  33  (2000),  no. 4, pp. 507--517.

\vskip 5pt

\noindent [13] McCrory, C.  {\it Zeeman's filtration of homology}.  Trans. Amer. Math. 
Soc. 250  (1979), pp. 147--166.

\vskip 5pt

\noindent [14] Milnor, J.  {\it Singular points of complex hypersurfaces}.
Annals of Mathematics Studies, No. 61 Princeton University Press, Princeton, N.J.; University
of Tokyo Press, Tokyo  1968 iii+122 pp.

\vskip 5pt

\noindent [15] Neumann, W. D.  {\it The fixed group of an automorphism of a word hyperbolic
group is rational}.
 Invent. Math.  110  (1992),  no. 1, pp. 147--150.

\vskip 5pt

\noindent [16]  Swarup, G. A.  {\it On the cut point conjecture}.
 Electron. Res. Announc. Amer. Math. Soc.  2  (1996),  no. 2, pp. 98--100 (electronic).

\vskip 5pt

\noindent [17]  Zeeman, E. C.  {\it Dihomology. III. A generalization of the Poincar\'e
duality for manifolds}.
 Proc. London Math. Soc. (3)  13  (1963), pp. 155--183.

\end{document}